\documentclass[authoryear,a4paper,10pt,onecolumn,preprint,3p]{elsarticle}

\usepackage{url}
\usepackage{amsmath}
\usepackage{booktabs}
\usepackage{multirow}
\usepackage{caption}
\usepackage[normalem]{ulem}
\useunder{\uline}{\ul}{}
\newcommand*{\rom}[1]{\expandafter\@slowromancap\romannumeral #1@}
\newcommand{\RN}[1]{%
	\textup{\uppercase\expandafter{\romannumeral#1}}%
}
\usepackage{booktabs}

\usepackage[colorlinks=true, urlcolor=blue, pdfborder={0 0 0}]{hyperref}

		\newcommand{\git}{\mathbin{
				\mathchoice{/\mkern-6mu/}
				{/\mkern-6mu/}
				{/\mkern-5mu/}
				{/\mkern-5mu/}}}

\usepackage{lmodern}
\usepackage{amsmath}
\usepackage{mathtools}
\usepackage[ruled, vlined, linesnumbered]{algorithm2e}

\begin{document}
	\begin{frontmatter}
		\title{A two-stage stochastic approach for the asset protection problem during escaped wildfires with uncertain timing of a wind change}
		\author[add1]{Iman Roozbeh\corref{cor1}}
		\ead{iman.roozbeh@rmit.edu.au}
		\author[add1]{John Hearne}
		\ead{john.hearne@rmit.edu.au}
		\author[add2]{Babak Abbasi}
		\ead{babak.abbasi@rmit.edu.au}
		\author[add1]{Melih Ozlen}
		\ead{melih.ozlen@rmit.edu.au}
		
		\cortext[cor1]{Corresponding author}
		\address[add1]{School of Science, RMIT University, Melbourne, Australia}
		
		\address[add2]{School of Business IT \& Logistics, RMIT University, Melbourne, Australia}
		
		\begin{abstract} 
		 
		  	Wildfires are natural disasters capable of damaging economies and communities. When wildfires become uncontrollable, Incident Manager Teams (IMT's) dispatch response vehicles to key assets to undertake protective tasks and so mitigate the risk to these assets.	
		  	In developing a deployment plan under severe time pressure, IMT's need to consider the special requirements of each asset, the resources (vehicles and their teams), as well as uncertainties associated with the wildfire. A common situation that arises in southern Australian wildfires is a wind change. There is a reliable forecast of a wind change but some uncertainty around the timing of that change.
		  	To assist IMT's to deal with this situation we develop a two-stage stochastic model to integrate such an uncertainty with the complexities of asset protection operations. This is the first time a mathematical model is proposed which considers uncertainty in the timing of a scenario change. The model is implemented for a case study that uses the context of the 2009 Black Saturday bushfires in Victoria. A new set of benchmark instances are generated using realistic wildfire attributes to test the computational tractability of our model and the results compared to a dynamic rerouting approach. The computations reveal that, compared with dynamic rerouting, the new model can generate better deployment plans. The model can achieve solutions in operational time for realistic-sized problems although for larger problems the sub-optimal rerouting algorithm would still need to be deployed.

		\end{abstract}
		
		

	\end{frontmatter}
	

	\section{Introduction}
	
	Numerous wildfires in the Northern Hemisphere during the summer of 2018 have again underscored their devastating effect on lives and property. The July fires which began in Attica, Greece, for example, claimed the lives of 99 people. Wildfires are also a serious hazard throughout most of Australia. Some of the worst wildfires have occurred in Victoria including some areas near Melbourne. The so-called "Black Saturday" fires of 2009 was the most destructive, claiming 174 lives and destroying 2030 homes and over 3500 structures \citep{whittaker2009victorian}. 
	
	Incident Management Teams (IMT's) are responsible for managing the response to a wildfire event. Community safety is the primary objective with protection of assets also very important. IMT's take command of operations to mitigate the hazard and must deploy the resources available to achieve it. This aspect of the IMT's role is extremely complex and decisions need to be made under severe pressure of time. In this paper, we consider the case where a wildfire has developed beyond the possibility of suppression and resources need to be deployed to reduce the risk to critical assets. Examples of such assets include electrical sub-stations, bridges, and schools. The loss of these assets all have consequences which could affect a community for weeks or even months. Activities to reduce the risk of damage might include cutting trees and clearing debris around buildings, removing gas bottles and flammable items such as wooden furniture, and hosing down structures. Determining the optimal deployment of resources to minimise losses is known as the Asset Protection Problem (APP) \citep{van2014mixed}.
	
	The APP considers a landscape comprising assets at various locations connected by a network of roads some of which may be unsealed. Each asset has a set of particular requirements. For example, some assets are only accessible by 4x4 vehicles, and some assets have no nearby access to reticulated water and require a tanker for hosing down. The resources available to deploy the mitigating actions comprise different types of vehicles. Each vehicle type has a set of capabilities. To each asset IMT's must deploy vehicles with the appropriate capabilities to meet the requirements of that asset. For example, an asset with no reticulated water at the end of an unsealed road would need to be serviced by a tanker (rather than a pumper) with 4x4 capability. The problem is further complicated by the wildfire spreading across the landscape. This means that each asset must be serviced within a specific time-window. Obviously firefighters must depart before the arrival of the uncontrollable wildfire. On the other hand servicing an asset too early might lead to the drying out of dampened assets or the accumulation of further debris.
	
	The APP was first addressed by \cite{van2014mixed}. They formulated a mixed integer programming model to solve a deterministic problem where time windows were determined by fire spreading in a fixed direction with constant velocity. While this work served to highlight the utility of such an approach for the APP, solutions using commercial solvers could not solve realistic problems in operational time. Later \cite{roozbeh2018adaptive} developed a heuristic for this same problem that achieves good solutions within minutes and hence a valuable tool for IMT's. While this has proved useful there is a particular circumstance associated with wildfires as witnessed in south-eastern Australia which need addressing. In the "Black Saturday" and other wildfires in Victoria a change in wind direction was a major contributor to the devastation caused \citep{cruz2012anatomy, bradstock2012wildfires}. Hot, dry winds typically come from the northwest and are followed by a southwest wind change (\cite{cruz2012anatomy}). A long fire flank is then transformed into a large fire front \citep{mccaw2013managing}. The updraft effect on the heated air by the cooler southwesterly also leads to an increase in the spread rate of the fire. Meteorologist are able to forecast that a change in wind direction is coming but there is a window of uncertainty about the exact time that this will occur. In this paper we extend previous work on the APP to deal specifically with this aspect of a wind change when the time of change is uncertain.
	
	Although there are many studies on the dispatch of resources for direct fire suppression (see \cite{duff2015operational, o2016getting, o2017empirical}), there are still very few studies on the APP\citep{van2014mixed, roozbeh2018adaptive} and other defensive tasks (\cite{donovan2003integer, minas2012review, van2017dynamic}). However, the APP can be regarded as a special case of the vehicle routing problem (VRP) as discussed below.
	
	The vehicle routing problem (VRP) and its variants has been intensively studied and future research directions are set in the literature (\cite{cordeau2007vehicle, toth2014vehicle, laporte2009fifty}). \cite{Archetti2014273} provide a survey on the VRP with profit in which the orienteering problem (OP) was considered as the basic problem of this class. The OP is a combination of the travelling salesman problem and the knapsack problem where travelling to all nodes is often not feasible due to a time constraint (\cite{ vansteenwegen2011orienteering, gunawan2016orienteering, roozbeh2016heuristic}). The APP is analogous to the OP. Constrained by time, both problems involve choosing a subset of all nodes to visit to maximise an objective. In the case of the OP a reward is on offer at each node and the objective is to maximise the total rewards collected. The objective of the APP is to maximise some weighted total of the number of sites serviced. Sites are often weighted by a measure of the economic consequence of the facility being destroyed or damaged. Unlike the OP, the APP may require synchronous visits to some nodes (assets) by more than one vehicle type. Dynamic fire fronts further impose time window constraints on the service time of each node as well as on the accessibility of roads. 
	
	Recently, research interest on stochastic variants of routing problems has increased significantly. Advances in technology have enabled larger, more complex problems to be solved to support decision makers. For further details, we refer the interested readers to the surveys by \cite{gendreau201650th} and \cite{ritzinger2016survey}. When dealing with stochastic problems, a wide range of approaches can be utilised, e.g. discrete event simulation or robust optimisation (\cite{hoyos2015or}). One of the most frequently used techniques that has attracted substantial attention in stochastic VRP and disaster management problems is the two-stage stochastic programming (\cite{falasca2011two, grass2016two, krasko2017two}). In brief, the two-stage stochastic programming make decisions in the first-stage prior to realisation of uncertain events, while taking prospective second-stage decisions into account. 
	There are various stochastic characteristics of the orienteering problem that have been investigated. \cite{ilhan2008orienteering} introduced the OP with stochastic profit for the first time. Other uncertainties, such as stochastic travel and service times (\cite{papapanagiotou2014objective}), stochastic time-dependent travel times (\cite{varakantham2013optimization}) and stochastic waiting time \cite{zhang2014priori} have been studied. While to the best of our knowledge, no formulations have been reported in the literature of a two-stage stochastic programming model with an uncertain time of change, i.e. staging time. Staging time is the moment when transition from one stage to another occurs. This is a situation that is commonly faced in emergency and logistic problems.
	
	This paper handles the uncertainty in the time of a change by developing a novel two-stage stochastic programming model. The proposed optimisation model incorporates the fact that the advancing fire directions will change but there is uncertainty about the timing of this change. 
	
	The first-stage plans operations according to the time that the direction of the wildfire spread might change.	
	The number and type of vehicles required to service each site is a particular requirement of the asset at that site. If more than one vehicle is required they must visit synchronously and within the time windows imposed by the advancing fire front.
	Deployed vehicles, in general, are unable to service all assets and so assets must be prioritised according to a predetermined economic value on the community's asset register. The model yields a deployment plan to maximise the expected value of the 'at risk' assets protected in all scenarios assuming that each asset serviced is protected. The first stage operations consider assets that are likely to be affected in different scenarios of the second stage. in other words, due to the uncertain time of wind change, the first stage tasks are being planned to improve the expected value of the protected assets under different stage transition scenarios.

	The remainder of this paper is structured as follows: we describe and develop a two-stage stochastic programming model for the asset protection problem in Section \ref{probdef}. Section \ref{compu}. describes a case-study and computational results are presented afterwards. In Section \ref{extension}. an extension to the model is identified and explained. Finally, the concluding remarks and future research directions are presented in Section \ref{conclusion}.

	\section{Problem description and model formulation} \label{probdes}
	Our two-stage stochastic approach for the APP encompasses multiple characteristics of the problem, such as a heterogeneous fleet of trucks, multiple scenarios, uncertain time of change and locations, time windows defined by the fire front, and synchronous service requirements. The problem is formulated in a generic manner to facilitate its application to analogous problems.

	\subsection{Sets, parameters and decision variables}

	The mixed integer linear programming formulation uses the following notation.
	
	\textbf{Indices and sets}

	$Q\;\;\;\;\;\;\;\;\;\;$ 			set of vehicle types
	
	$A\;\;\;\;\;\;\;\;\;\;$ 			set of all arcs 
	
	$\delta^+_q(i)\;\;\;\;\,$		set of feasible arcs $(i,j)$ that can be traversed from $i$, $q \in Q$
	
	$\delta^-_q(j)\;\;\;\;$		set of feasible arcs $(i,j)$ that can be traversed to $j$, $q \in Q$
	
	$N\;\;\;\;\;\;\;\;\;\,$ 	    	set of all assets
	
	
	$\Xi\;\;\;\;\;\;\;\;\;\;\,$ 	    	set of all scenarios, $\xi \in \Xi$

	\textbf{Parameters}
	
	$a_i\;\;\;\;\;\;\;\;\;\;\,$ 			service duration associated with location $i$

	$c_i^f\;\;\;\;\;\;\;\;\;\;$ 			latest time that protection activities may commence at stage $f$ (first stage)
	
	$c_i^s(\xi)\;\;\;\;\;\;$ 			latest time that protection activities may commence at stage $s$ (second stage) in scenario $\xi \in \Xi$

	$P(\xi)\;\;\;\;\;\;\,$ the probability that scenario $\xi \in \Xi$ occurs
	
	$n\;\;\;\;\;\;\;\;\;\;\;\,$ 				number of assets in the graph representation of the problem
	
	$o_i^f\;\;\;\;\;\;\;\;\;\;$ 			earliest time that protection activities may commence at location $i$ in stage $f$
	
	$o_i^s(\xi)\;\;\;\;\;\;$ 			earliest time that protection activities may commence at location $i$ in stage $s$, scenario $\xi \in \Xi$
	
	$\kappa_q\;\;\;\;\;\;\;\;\;\;$ 				number of vehicles of type $q \in Q$
	
	$r_i\;\;\;\;\;\;\;\;\;\;\;$ 				protection requirement vector for asset $i$, e.g. $r_i = <2, 1, 0>$
	
	$start_{p}\;\;\;\;$ 			1 if vehicle $p$ is at depot, 0 otherwise
	
	$t_{ijq}\;\;\;\;\;\;\;\;\,$ 			travel time from location $i$ to location $j$ by vehicle $q \in Q$
	
	$\nu_i\;\;\;\;\;\;\;\;\;\;\,\,$ 				value of asset $i$
	
	$ST\;\;\;\;\;\;\;\;\;$ 				staging time
	
	$T_{max}\;\;\;\;\;\;$ 				latest time allowed to start a protective task 
	
	$0, N+1\;$ 		represents the start and final locations (may be the same location)

	\textbf{Variables}
	
	$w_{j}\;\;\;\;\;\;\;\;\;\;$ 				1 if asset $j$ is visited just before stage transition occurs, 0 otherwise

	\textbf{First stage decision variables}
	
	$X_{ijq}^f\;\;\;\;\;\;\;$ 			number of vehicles of type $q \in Q$ travelling from location $i$ to location $j$ at stage $f$
	
	$z_{ijq}^f\;\;\;\;\;\;\;\;$ 				1 if arc $(i,j)$ is traversed by vehicle $q \in Q$, 0 otherwise
	
	$S_i^f\;\;\;\;\;\;\;\;\;\;$ 				time at which service commences in location $i$ at stage $f$
	
	$Y_i^f\;\;\;\;\;\;\;\;\;$ 				1 if location $i$ is serviced at stage $f$, 0 otherwise

	\textbf{Second stage decision variables}

	$z_{ijq}^s(\xi)\;\;\;\,$ 				1 if arc $(i,j)$ is traversed by vehicle $q \in Q$ at stage $s$ in scenario $\xi \in \Xi$, 0 otherwise
	
	$Y_i^s(\xi)\;\;\;\;\:$ 				1 if location $i$ is serviced at stage $s$ in scenario $\xi$, 0 otherwise

	$S_i^s(\xi)\;\;\;\;\;\;$ 				time at which service commences in location $i$ at stage $s$ in scenario $\xi \in \Xi$

	$X_{ijq}^s(\xi)\;\;\;$ 			number of vehicles of type $q \in Q$ travelling from location $i$ to location $j$ at stage $s$ in scenario $\xi $
	
	$Y_i^s(\xi)\;\;\;\;\;$ 				1 if location $i$ is serviced at stage $s$ in scenario $\xi \in \Xi$, 0 otherwise
	
	\subsection{Mathematical model} \label{sec:formulation}
	
	 Consider a set of $n$ assets with depots $0$ and $N+1$. The start and end depots may be the same. Each asset $i$ is associated with a value $\nu_i$ and a vector of service requirements $R_i=<r_{i1},r_{i2},\dots,r_{iq}>$. $R_i$ indicates the number of vehicles of each type needed to accomplish a task cooperatively to service an asset.
	Variables $Y_i$ is binary and equal 1 if an asset is protected and 0 otherwise. 
	The decision variable $z_{ijq}$ takes value 1 if arc $(i,j)$ traversed by vehicle $q$ while variables $X_{ijq}$ and $t_{ijq}$ represent number of vehicles travel along arc $(i,j)$ and the travel time.
	Superscripts $f$ and $s$ indicate the stage and $P(\xi)$ is the probability of scenario $\xi$ at the second stage. Every asset should be visited within the associated time window $[o_i,c_i]$ depending on the scenario and the stage. The set of all arcs is defined by $A$. The set $\delta^+_q(i)$ is the set of all feasible arcs $(i,j)$ for each asset $i$ such that both assets $i$ and $j$ need to be visited by vehicle type $q\in Q$ and can be reached within their time windows. Similarly, $\delta^-_q(j)$ is the set of all feasible incoming arcs to asset $j$. Parameter $ST$ represents the staging time where a change in the problem conditions occurs. A stage transition may occur at any location on an arc or at an asset where a vehicle is located when the staging time occurs. $w_i$ is assigned the value 1 if $i$ is the last asset on a route to be visited before the commencement of stage two. This enables us to maintain information through the stage transition. 
	The superscripts $f$ and $s$ indicate variables for the first and second stage, respectively.
	
	With the notations and explanations above, the problem is formulated as follows:

	
	\begin{equation}\label{1}	 	
	Max \;\;\; \sum_{i\in N } \nu_{i} Y_{i}^f + 
	\sum_{\xi \in \Xi } P(\xi)(\sum_{i\in N }\nu_{i} Y_{i}^s{(\xi)})
	\end{equation}
	
	\begin{equation}\label{2}	 
	s.t.:\sum_{j \in \delta^+_q(0)} X_{0jq}^f
	=
	\sum_{i \in \delta^-_q(N+1)} X_{i(N+1)q}^f
	+
	\sum_{i \in \delta^-_q(N+1)} X_{i(N+1)q}^s(\xi)
	, \qquad \forall q\in Q, \forall \xi  \in \Xi;
	\end{equation}
	
	\begin{equation}\label{3}	 
	\sum_{i \in \delta^-_q(j)}X_{ijq}^f + \sum_{i \in \delta^-_q(j)}X_{ijq}^s(\xi) = \sum_{k \in \delta^+_q(j)}X_{jkq}^f + \sum_{k \in \delta^+_q(j)}X_{jkq}^s(\xi)   , \qquad \forall q\in Q,\forall j\in N,\forall \xi  \in \Xi;
	\end{equation}
	
	\begin{equation}\label{4}	 
	\sum_{q\in Q}\sum_{k \in \delta^+_q(j)}X_{jkq}^f -\sum_{q\in Q} \sum_{i \in \delta^-_q(j)}X_{ijq}^f \ge -M \times w_j  , \qquad  \forall j\in N;
	\end{equation}
	
	\begin{equation}\label{5}	 
	\sum_{q\in Q}\sum_{k \in \delta^+_q(j)}X_{jkq}^f - \sum_{q\in Q}\sum_{i \in \delta^-_q(j)}X_{ijq}^f \le -w_j  , \qquad \forall j\in N;
	\end{equation}
	
		\begin{equation}\label{6}	 
	S_j^{f} + a_j - ST\le M (1-z_{ijq}^f)  , \qquad \forall q\in Q, \forall (i,j)\in A;
	\end{equation}
	
		\begin{equation}\label{7}	 
	S_j^{s}(\xi)+a_j - ST\ge M (z_{ijq}^s(\xi)-1)  , \qquad \forall q\in Q, \forall (i,j)\in A, \forall \xi  \in \Xi;
	\end{equation}
	
		\begin{equation}\label{8}	 
	\sum_{j\in\delta_q^+(k)}X_{0jq}^f \le start_{q} , \qquad \forall q\in Q;
	\end{equation}
	
		\begin{equation}\label{9}	 
	\sum_{i\in\delta_q^-(j)}X_{ijq}^f = r_{jq} Y_j^f, \qquad \forall j\in N,  \forall q\in Q;
	\end{equation} 	
	
	\begin{equation}\label{10}	 
	\sum_{i\in\delta_q^-(j)}X_{ijq}^s(\xi) = r_{jq} Y_j^s(\xi), \qquad \forall j\in N,  \forall q\in Q, \forall \xi  \in \Xi;
	\end{equation} 
	
	\begin{equation}\label{11}	 
	Y_{j}^f + Y_{j}^s(\xi) \le 1, \qquad \forall j\in N,  \forall q\in Q, \forall \xi  \in \Xi;
	\end{equation} 
	
	\begin{equation}\label{12}	 
	X_{ijq}^f \le \kappa_{q} z_{ijq}^f, \qquad \forall (i,j)\in A,\forall  q\in Q; 
	\end{equation} 	 	
	
		\begin{equation}\label{13}	 
	X_{ijq}^s(\xi) \le \kappa_{q} z_{ijq}^s(\xi), \qquad \forall (i,j)\in A,\forall  q\in Q, \xi  \in \Xi;
	\end{equation} 
	
	\begin{equation}\label{14}	 
	S_{i}^{f}+t_{ijq}+a_i-S_j^{f} \le M(1-z_{ijq}^f), \qquad \forall (i,j)\in A, \forall q\in Q;
	\end{equation}
	
	\begin{equation}\label{15}	 
	S_{i}^{s}(\xi)+t_{ijq}+a_i-S_j^{s}(\xi) \le M(1-z_{ijq}^s(\xi)), \qquad \forall (i,j)\in A, \forall q\in Q, \forall \xi  \in \Xi;
	\end{equation} 	
	
	\begin{equation}\label{16}	 
	S_{i}^{f} - S_{i}^{s}(\xi)\le  M  (1 - w_i), \qquad \forall (i)\in N, \forall \xi  \in \Xi;
	\end{equation} 
	
	\begin{equation}\label{17}	 
	S_{i}^{f} - S_{i}^{s}(\xi)\ge  M  ( w_i - 1), \qquad \forall (i)\in N, \forall \xi  \in \Xi;
	\end{equation} 
	
	\begin{equation}\label{18}	 
	o_{j}^f \le S_j^{f}\le c_{j}^f, \qquad \forall j\in N;
	\end{equation} 	
	
	\begin{equation}\label{19}	 
	o_{j}^s(\xi) \le S_j^{s}(\xi)\le  c_{j}^s(\xi), \qquad \forall j\in N, \forall \xi \in \Xi;
	\end{equation} 	
	
	\begin{equation}\label{20}	 
	Y^f_i, Y_{i}^s{(\xi)}\in \{0,1\}, \qquad \forall i\in N, \forall \xi  \in \Xi;
	\end{equation}
	
	\begin{equation}\label{21}	 
	z^f_{ijq}, z_{ijq}^s{(\xi)}\in \{0,1\}, \qquad \forall (i,j)\in A, \forall q \in Q, \forall \xi \in \Xi;
	\end{equation}
	
	\begin{equation}\label{22}	 
	w_i\in \{0,1\}, \qquad \forall i\in N;
	\end{equation}
	
	\begin{equation}\label{23}	 
	X^f_{iiq}, X_{ijq}^s{(\xi)} \in \{0,1,\dots,\kappa_q\}, \qquad \forall (i,j)\in A, \forall q \in Q, \forall \xi \in \Xi;
	\end{equation}
	
	\begin{equation}\label{24}	 
	S_i^f, S_{i}^s{(\xi)} \in [0,T_{max}], \qquad \forall i\in N, \forall \xi \in \Xi.
	\end{equation}
	The objective function (\ref{1}) maximises the expected value of serviced assets.
	Constraints (\ref{2}) require that all vehicles leaving the starting depot must reach the final depot.
	Constraints (\ref{3}) enforces the conservation of flow at each asset.
	Constraints (\ref{4} and \ref{5}) indicate the staging location. If a vehicle depart a node and cannot arrive to the destination by the staging time, the departing asset will be denoted as staging location. This assists in assigning the right values to $S_j^{s}(\xi)$ in constraints (\ref{16} and \ref{17}) in order to start trips from staging locations in various scenarios at stage two.
	Constraints (\ref{6} and \ref{7}) impose the staging time. If an arc is traversed in stage $f$, $z_{ijq}^s(\xi)$ for all scenarios must be zero. On the other hand, if $z_{ijq}^f$ is zero, one of the $z_{ijq}^s(\xi)$ could be one, but not necessarily.
	Constraints (\ref{8}) indicate that the number of vehicles departing from a depot may not exceed the number of vehicles stationed at the depot. 
	Constraints (\ref{9} and \ref{10}) guarantee an asset is serviced only if its service requirements are fulfilled by the right combination of incoming vehicles.
	Constraints (\ref{11}) ensure each asset will be visited only once in both stages.
	Constraints (\ref{12} and \ref{13}) make sure that vehicles travelling through an arc never exceed the number available $\kappa_q$.
	Constraints (\ref{14} and \ref{15}) ensure that an asset may only be visited if the service requirements of the previous location has been satisfied and there is enough time to reach the next asset.
	Equations (\ref{16} and \ref{17}) guarantee that proper values are assigned to $S_j^{s}(\xi)$ through a transition from one stage to another.
	Terms (\ref{18}) and (\ref{19}) ensure that the time window constraints are not	violated.
	Constraints (\ref{20}-\ref{24}) impose non-negative, binary and integer restrictions on the variables.

	\subsection{Fire spread} \label{fire}
	
	The time windows during which an asset must be serviced will depend on the progress of a wildfire as it spreads across the landscape. There are several models available for simulating the advance of a wildfire. Reviews can be found in \cite{scott2005standard,sullivan2009wildland,johnston2008efficient} while \cite{Petrasova2018firemodel} is an example of recent developments. Wildfire behaviour is dependent on factors such as temperature, humidity, fuel load, etc. To illustrate and then test our model we need to simulate the pattern of spread of a wildfire across an artificial landscape. From this we can determine the time windows for servicing each asset in the landscape. Due to the unique circumstances of each wildfire only some general features of wildfire spread are used to generate our benchmark set.
	With a constant wind fire-fronts mostly progress in an elliptical shape (\cite{anderson1982modelling}). Thus we use the general equation of an ellipse to determine the time that fire impacts an asset. 

\begin{figure}[h!]
		\centering
		\includegraphics[width=1\textwidth,height=0.5\textheight]{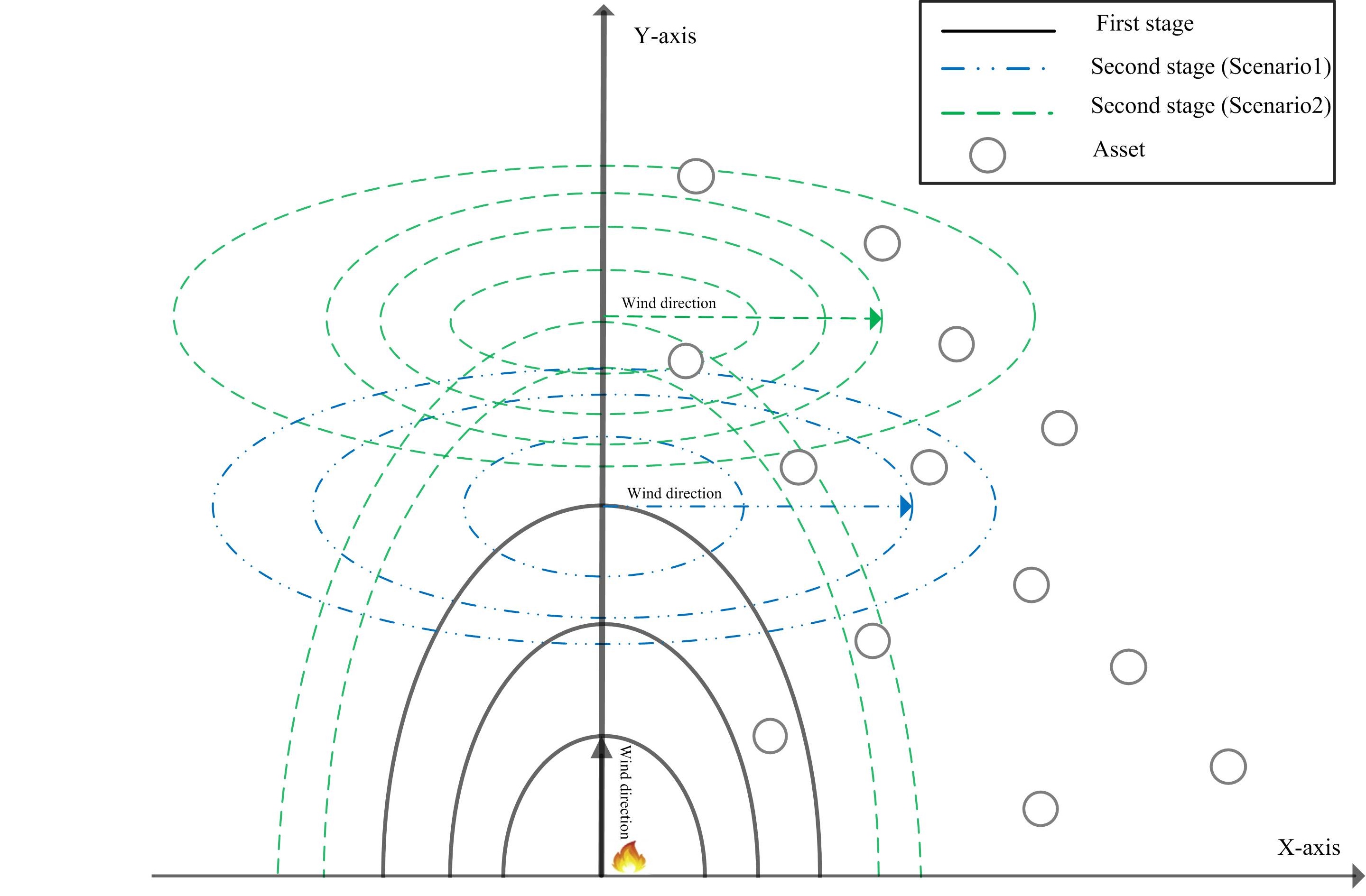}
		\caption{Fire front progressing through a landscape. The initial north wind (Y-direction) is expected to change to a westerly but there is uncertainty as to the timing of the change. Two scenarios are shown representing the times at which the wind change occurs. In each scenario a new set of assets are at risk.}
		\label{fire-spread}
	\end{figure}

	The problem we want to investigate is the frequent situation where it is known that a wind change will occur but in stage 1 there is uncertainty over its timing. As shown in Figure \ref{fire-spread} the timing of the wind change is represented by two scenarios. In the first scenario the wind change occurs at the staging time, $(ST)$, and an hour later in scenario 2. The change in wind direction and its timing has a significant impact on the fire front. This in turn affects which assets will be in the path of the wildfire.
 
	\newpage
	\subsection{An illustrative example}
\begin{figure}[h!]
		\centering
		\includegraphics[width=0.7\textwidth,height=0.3\textheight]{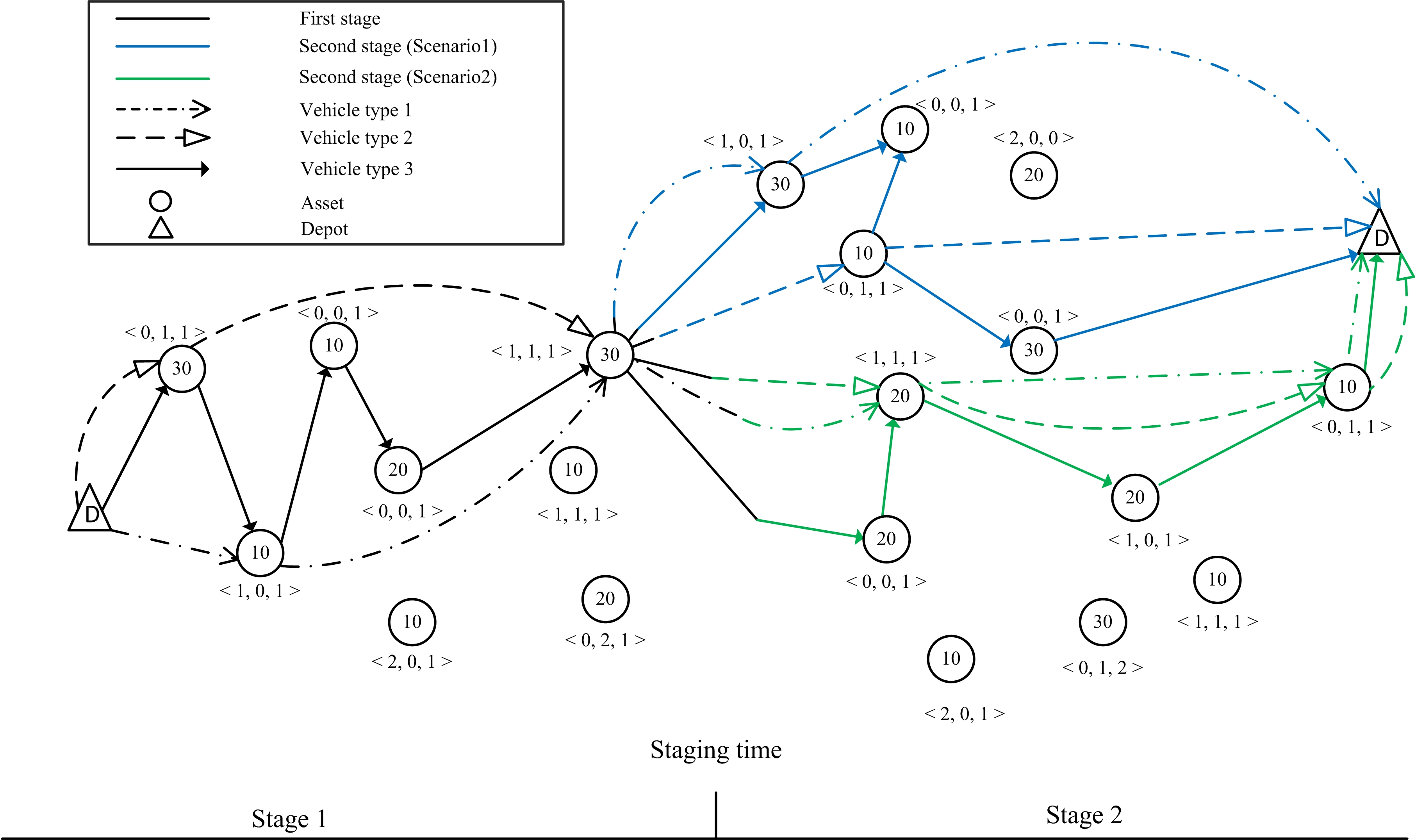}
		\caption{Illustrative APP problem with two wind change scenarios. The area impacted in each scenario is shown in different colours. Asset values are shown within the circles with the requirements shown below them. Vehicle types 1, 2 and 3 have capability vectors of $<1,0,0>$, $<0,1,0>$ and $<0,0,1>$. To the extent possible they must be assigned to the assets so that their capabilities match the asset's requirements. The optimal assignment of the three vehicles is shown.}
		\label{probdef}
	\end{figure}

	A simple demonstration of our two-stage optimisation model for the APP problem is sketched in Figure \ref{probdef}. The time windows were generated using the method described in Section \ref{fire} above. The service requirements of each asset were chosen to highlight some of the key features of the APP problem. Due to resource limitations and time windows not all assets can be serviced. The optimisation assigned vehicles to maximise the expected value of assets serviced. The assets actually serviced also depends on which scenario is realised.

	As shown in Figure \ref{probdef}, twenty community assets are at risk. Each vehicle and asset has particular capabilities and service requirements that need to be matched. For example, an asset in stage 1 which has requirements of $<2,1,1>$ must be serviced by two type 1 vehicles and one of the other two vehicle types. The service must start cooperatively and simultaneously within the time window $[o_i^f,c_i^f]$ for a duration of $a_i$. In the Figure \ref{probdef}, the probability of occurrence for each scenario is known and routing in the first stage is planned to gain the highest benefit from the expected value of the objective function in the second stage. 

	\newpage
	\section{Experiments} \label{compu}
	We carried out extensive computational experiments to test our model. Initially, we solve a case study using the context of the February 2009, Black Saturday bushfires in Australia. Then, we perform further tests on a large set of generated benchmark instances using realistic parameters \footnote{All the benchmark instances are available at www.sites.google.com/site/imanrzbh/datasets} and compare them against dynamic rerouting approach. All computational experiments are implemented on a desktop computer equipped with Intel Core i5 (3.2GHz) and 8.0 GB of RAM, where MILP models are solved by the commercial solver, CPLEX 12.8, coded in Python 3.6. Computational times are measured in elapsed time with a time limit of an hour.
	
	\subsection{Case study - Murrindindi Mill fire Black Saturday} \label{case}
	
\begin{figure}[h!]
		\centering
		\includegraphics[width=1\textwidth,height=0.5\textheight]{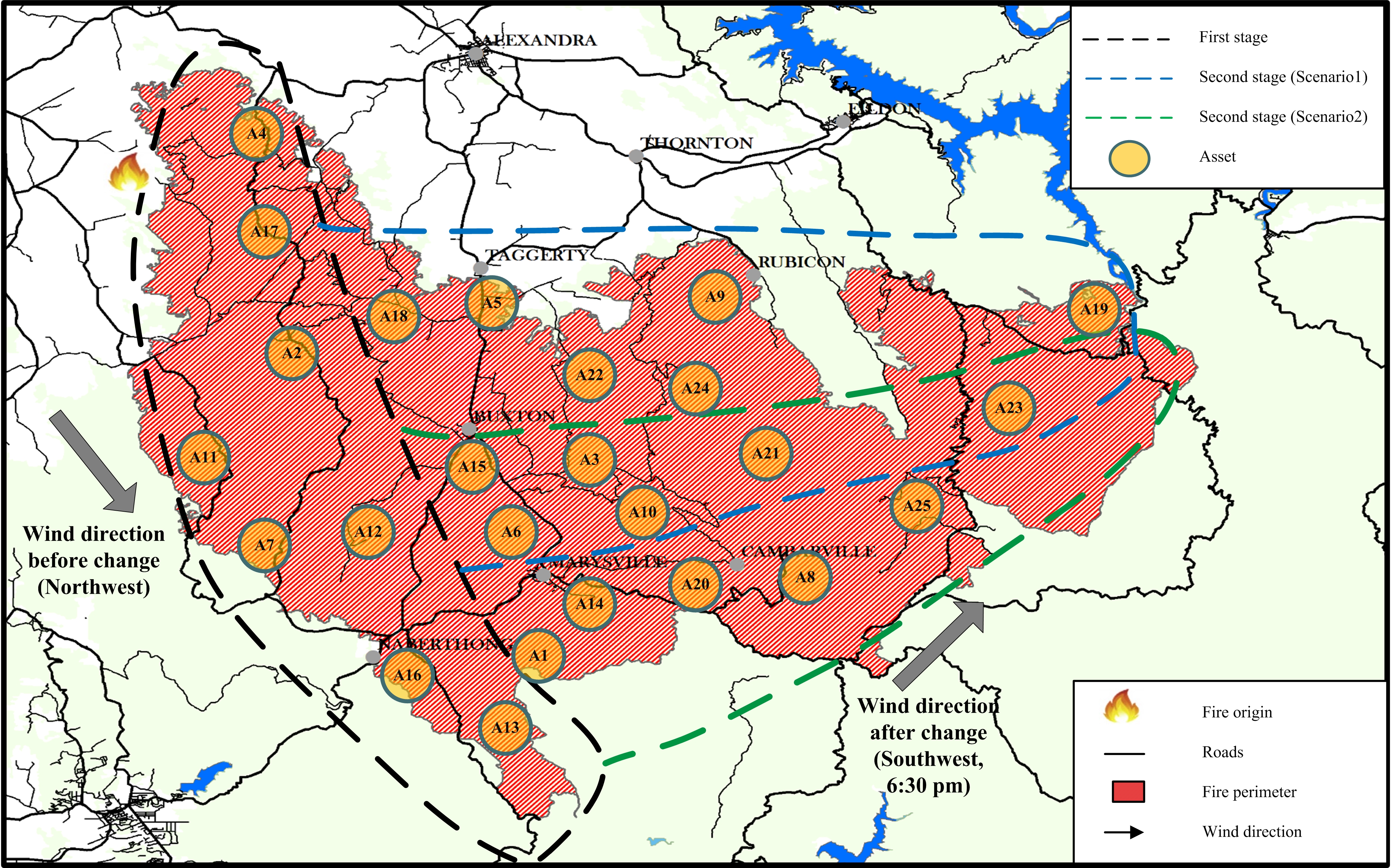}
		\caption{Case study region - Murrindindi Shire, Victoria, Australia.}
		\label{blacksat}
	\end{figure}

	The shire of Murrindindi is $3,889 \; km^2$ in extent with a population of $13,732$ (2016 census) located in the north-eastern part of Victoria, Australia. About $46\%$ of the total land area of the municipality is forest (1788 km2), and a large proportion of this land is mountainous (\cite{shahparvari2016enhancing}). On 7 February 2009, a series of bushfires known as Black Saturday raged through the shire (see Figure \ref{blacksat}). The fire swept the 50 km distance between Saw Mill in Wilhelmina Falls Road and Narbethong in about 90 minutes. Following a wind change the long narrow fire-front become a wide fire-front that burned through a number of townships with tragic consequences (\cite{cruz2012anatomy}). Five of the bushfires on Black Saturday claimed people's lives. The second highest number of deaths resulted from the Murrindidi bushfires (40 people), which reveals the importance of the area under study. (\cite{whittaker2009victorian}).

	A set of 25 assets are identified in the area as being impacted (see \ref{appendixB}). The Yea country fire authority is defined as the starting point of the operations (depot). Distance matrices were obtained using the Google Maps API (\cite{GoogleMapsAPI}). Time windows are set approximately based on the Victorian Bushfires Royal Commission report (see \cite{teague2010final}). The number of vehicles are considered proportional to the problem size in order to cover a significant amount of at risk assets $\{\kappa_1 = 5,\,\kappa_2 = 3,\,\kappa_3 = 2\}$. Vehicle velocity ($V=60\,km.h^{-1}$) and probability of scenarios ($P(1)=\;70\%,\, P(2)=\;30\%$) are other parameters involved in the problem. The last piece of the following case study is the staging time ($ST = 6:30\,pm$), the time that wind direction changes and new assets will be impacted. 
	
	We used CPLEX to solve the case study discussed above using the two-stage stochastic programming model formulated in Section \ref{sec:formulation} \footnote{Asset information and detailed solutions are available at www.sites.google.com/site/imanrzbh/datasets}. The optimal solution showed that only two assets, $11$ and $13$, under threat during stage 1 could not be serviced. A total of ten vehicles finished their protective activities at assets $4$, $12$, $16$ and $17$ before immediately heading to the assets at risk in the next stage. 
	
	It can be seen in Figure \ref{blacksat} that some assets, $3$, $6$, $10$, $15$, $21$ and $23$, are impacted by the wildfire in both scenarios. The lack of uncertainty with assets appearing in both scenarios leads to a large proportion of them being serviced. This is true in this case for $4$ out of $6$ assets. Of course, the optimal solution will also be influenced by the value of an asset and the resources required. Overall, the optimal solution provides a strategy for servicing assets to the value of $86\%,\,71\%$ and $83\%$ of the total asset value in stage 1, and scenarios 1 ans 2, respectively.
	
	\subsection{Benchmark instances}\label{bench}
	To further evaluate the model, a set of benchmark instances are generated by randomly distributing assets over an $80\times80$ grid. Six sets of problems with asset numbers from 50 to 75 are solved. Each of these problems is solved for two cases representing different sets of vehicles as given in Table \ref{result}. We use the elliptical fire spread data explained in Section \ref{fire} to set time windows for each asset. Similar to the events of Black Saturday we consider a wind change. Two scenarios are possible for the timing of the wind change. In scenario 1 the wind change occurs at the staging time whereas it occurs two hours later in scenario 2. The realisation of which scenario will occur will only be known at the staging time but deployment decisions must be made at the start of stage 1.
    
    For representing a realistic situation where the impact of the wind change is significant, a different set of assets are impacted with each scenario. Therefore, there are four types of assets: (1) assets impacted at the first stage; (2) assets in either scenario 1 or scenario 2 but not both at the second stage; (3) assets that are affected in both scenarios, (4) and assets that are not at risk at all in either scenario. Fire velocities are set in a manner to simulate escaped wildfires. Travel time is calculated with vehicle velocity set at a conservative but realistic ($30 \: km.h^{-1}$). 
    To investigate the impact of vehicle numbers on the percentage of the assets protected, we performed experiments with two different sets of vehicles 
    $Set1=\{\kappa_1=3,\:\kappa_2=2,\:\kappa_3=2\}$ and $Set2=\{\kappa_1=4,\:\kappa_2=3,\:\kappa_3=2\}$.
	The key parameters for generating benchmark instances are mostly inspired by real life scenarios, and are listed in Table \ref{parameters}.

	To evaluate the performance of our model, we also solve each benchmark instance using 'dynamic rerouting' and compare the solutions of the two approaches.

	\subsubsection{Dynamic rerouting approach}\label{dy}
	
	The dynamic rerouting of vehicles reallocate resources in the event of a change in conditions. The initial deployment of vehicles is based on a deterministic MILP (see \ref{appendixA1}) assuming no changes will occur. If any change in conditions or disruptions occur during operations, new parameters (e.g. time windows) are updated accordingly. Then, the problem is solved to optimality from where the disruption took place. The pseudo-code of the dynamic rerouting asset protection problem is described in Algorithm \ref{ALG1}.
	
			\makeatletter
		\def\BState{\State\hskip-\ALG@thistlm}
		\makeatother
		\begin{algorithm} 
		    \caption{Pseudocode for the dynamic rerouting approach}\label{ALG1}
			\DontPrintSemicolon
			\KwIn{\\ staging time ($ST$), set of at risk assets in stage 1 ($N^f$), set of at risk assets in stage 2 under scenario \# 1 ($N_1^s$), set of at risk assets in stage 2 under scenario \#2 ($N_2^s$), probability of scenario \#1 ($P(1)$), probability of scenario \#2 ($P(2)$), collected rewards from assets in set $N^f$ ($\nu^f$), collected rewards from assets in set $N^s_1$ ($\nu^s_1$), collected rewards from assets in set $N^s_2$ ($\nu^s_2$), set of locations where resources are at time ST ($D^{st}$), deterministic asset protection model ($\Omega$), set of all assets ($N$), set of visited assets ($\Upsilon$).
		}
			\KwOut{ Expected value of the collected rewards ($E(\nu)$).}
			\SetKwBlock{Begin}{Function}{end function}
			\Begin($Dynamic \; rerouting$)
			{				
			\uIf{$(P(1)>P(2))$}{
				
				$N_a$ $\gets $ $ N^f \cup N^s_1$\;
				$N_b$ $\gets $ $  N^s_2$
			}
	\Else {
		\uIf {($(P(2)>P(1))$)}{
		$N_a$ $\gets $ $ N^f \cup N^s_2$\;
			$N_b$ $\gets $ $  N^s_1$
		}	
	}
			Solve $\Omega$ for $assets\; \in \; N_a$ $\qquad \git return\: \Upsilon$\;
			
			Define $D^{st}$
			
			Solve $\Omega$ for $assets\; \in \; N_b$ starting from $D^{st}$ $\qquad \git update\: \Upsilon$\;
			
			$E(\nu) \gets 0$
			
				\ForAll { ($i\in N $)}{
					\uIf{$(i \in \Upsilon \;\&\&\; i \in N^f)$}{
					$E(\nu) \gets E(\nu)+ \nu_i$	
					}
					\uIf{$(i \in \Upsilon \;\&\&\; i \in N^s_1)$}{
						$E(\nu) \gets E(\nu)+ P(1) \times \nu_i$	
					}
					\uIf{$(i \in \Upsilon \;\&\&\; i \in N^s_2)$}{
						$E(\nu) \gets E(\nu)+ P(2) \times \nu_i$	
					}
			}
				\Return {$E(\nu)$}
			}
			\end{algorithm}
		
	
	Algorithm \ref{ALG1}. returns the expected value of the protected assets using the dynamic rerouting approach. The algorithm defines a scenario that has the highest probability (lines 2-8). Then, constructs a set of nodes including the first stage assets and those belonging to the scenario with the highest probability of occurring. After solving the deterministic MILP, the model determines the set of protected assets. By solving the initial MILP at line 9 the algorithm also determines locations where the resources would be at when the disruptions occur. Algorithm \ref{ALG1}. solves the MILP model again for assets belonging to other scenarios, while resources start their routes from where they were at the staging time. In lines (13-19) the algorithm calculates the expected value of the assets protected by multiplying the total asset value at each scenario by the probability of each scenario.
	The same benchmark instances are used and solved according to the above procedure. Results are reported in Table \ref{result}.

	\begin{table}[h!]
		\centering
		\caption{Parameters used in the benchmark instances}
		\label{parameters}
		\resizebox{1\columnwidth}{!}{%
			\begin{tabular}{@{}llc@{}}
				\toprule
				\multicolumn{1}{c}{\textbf{Parameter}} & \multicolumn{1}{c}{\textbf{Explanation}}                        & \textbf{Value} \\ \midrule
				$v_{x0}$                                    & Fire velocity along x in stage one.                             & 14             \\
				$v_{y0}$                                    & Fire velocity along y in stage one.                             & 16             \\
				$v_{x1}$                                    & Fire velocity along x in stage two, scenario one.               & 19             \\
				$v_{y1}$                                    & Fire velocity along y in stage two, scenario one.               & 17             \\
				$v_{x2}$                                    & Fire velocity along x in stage two, scenario two.               & 21             \\
				$v_{y2}$                                    & Fire velocity along y in stage two, scenario two.               & 19             \\
				$Delay$                                  & Delay between staging time and change of wind, in scenario two. & 2              \\
				$TW_1$                                    & Length of time window for assets impacted in stage one.          & 1              \\
				$TW_2$                                    & Length of time window for assets impacted in stage two.          & 1            \\
				$a$                                      & Service duration time.                                 & 0.5            \\
				$P(1)$                                     & Probability of scenario 1.                                      & 0.6            \\
				$P(2)$                                     & Probability of scenario 2.                                      & 0.4            \\
				$V$                                      & Travel speed of vehicles.                                       & 30             \\
				$ST$                                     & Staging time.                                                   & 4.5              \\ \bottomrule
			\end{tabular}
		}
	\end{table}
	
	\subsubsection{Numerical results} \label{results}
	
	Both the dynamic rerouting and the two-stage stochastic programming approaches are implemented to solve the benchmark instances. Note that we eliminate infeasible solutions due to the time window constraints, protection requirements and the impact time in a preprocessing step to sufficiently reduce the problem size for both approaches.	The result are presented in Table \ref{result}.

	\begin{table}[h!]
		\centering
		
		\caption{A summary of results for 50-55-60 number of assets. Vehicle numbers are defined at two levels: $Set1 = \{\kappa_1=3,\:\kappa_2=2,\:\kappa_3=2\}$ and $Set2 = \{\kappa_1=3,\:\kappa_2=3,\:\kappa_3=2\}$. Results are reported as percentage of the assets' values protected. Computation times are reported in seconds(sec).}
		\label{result}
		\resizebox{1\columnwidth}{!}{%
			\begin{tabular}{@{}llllllllllllll@{}}
				\cmidrule(l){1-14}
				
				\multirow{2}{*}{\#Vehicles} & \multirow{2}{*}{\#Assets} & \multicolumn{5}{l}{Two-stage stochastic programming}                                                                         &  & \multicolumn{5}{l}{Dynamic rerouting}                                                                                                                                                                                                                                                                                                                                                                        & \multicolumn{1}{c}{\multirow{2}{*}{\begin{tabular}[c]{@{}c@{}}GAP\\ (\%)\end{tabular}}} \\ \cmidrule(lr){3-8} \cmidrule(lr){9-13}
				&                           & \begin{tabular}[c]{@{}l@{}}Asset value\\ protected at\\ stage 1\\ (\%)\end{tabular} & \begin{tabular}[c]{@{}l@{}}Asset value\\ protected at\\ stage 2\\ scenario 1\\ (\%)\end{tabular} & \begin{tabular}[c]{@{}l@{}}Asset value\\ protected at\\ stage 2\\ scenario 2\\ (\%)\end{tabular} & \begin{tabular}[c]{@{}l@{}}Total value\\ of assets\\ protected\end{tabular} & \begin{tabular}[c]{@{}l@{}}Time\\ (sec)\end{tabular} &  & \begin{tabular}[c]{@{}l@{}}Asset value\\ protected at\\ stage 1\\ (\%)\end{tabular} & \begin{tabular}[c]{@{}l@{}}Asset value\\ protected at\\ stage 2\\ scenario 1\\ (\%)\end{tabular} & \begin{tabular}[c]{@{}l@{}}Asset value\\ protected at\\ stage 2\\ scenario 2\\ (\%)\end{tabular} & \begin{tabular}[c]{@{}l@{}}Total value\\ of assets\\ protected\end{tabular} & \begin{tabular}[c]{@{}l@{}}Time\\ (sec)\end{tabular} & \multicolumn{1}{c}{}                                                                    \\ \cmidrule(r){1-14} 
				\multirow{6}{*}{Set1}       & 50                        & 33.37                                                                               & 39.52                                                                                            & 33.29                                                                                            & 563.96                                                    & 39.41                                                &  & 33.74                                                                               & 39.52                                                                                            & 23.85                                                                                            & 532.52                                                    & 12.95                                                & 5.91                                                                                    \\
				& 55                        & 32.99                                                                               & 35.35                                                                                            & 28.51                                                                                            & 589.10                                                    & 70.12                                                &  & 33.37                                                                               & 35.48                                                                                            & 18.68                                                                                            & 544.41                                                    & 31.67                                                & 8.21                                                                                    \\
				& 60                        & 30.47                                                                               & 34.43                                                                                            & 26.39                                                                                            & 600.70                                                    & 202.11                                               &  & 30.64                                                                               & 34.43                                                                                            & 17.45                                                                                            & 550.44                                                    & 52.29                                                & 9.11                                                                                    \\
				& 65                        & 28.75                                                                               & 33.69                                                                                            & 25.34                                                                                            & 616.48                                                    & 359.14                                               &  & 29.04                                                                               & 34.02                                                                                            & 15.20                                                                                            & 560.96                                                    & 197.98                                               & 9.90                                                                                    \\
				& 70                        & 27.53                                                                               & 32.28                                                                                            & 23.35                                                                                            & 637.26                                                    & 555.48                                               &  & 27.77                                                                               & 33.10                                                                                            & 15.22                                                                                            & 576.33                                                    & 284.55                                               & 10.57                                                                                   \\
				& 75                        & 25.56                                                                               & 32.01                                                                                            & 23.23                                                                                            & 646.21                                                    & 1657.15                                              &  & 25.81                                                                               & 32.48                                                                                            & 15.29                                                                                            & 578.39                                                    & 617.56                                               & 11.73                                                                                   \\  \cmidrule(l){1-14} 
				\multirow{6}{*}{Set2}       & 50                        & 42.28                                                                               & 44.89                                                                                            & 34.32                                                                                            & 683.66                                                    & 178.49                                               &  & 42.28                                                                               & 45.44                                                                                            & 26.59                                                                                            & 639.07                                                    & 41.93                                                & 6.98                                                                                    \\
				& 55                        & 40.71                                                                               & 41.57                                                                                            & 33.76                                                                                            & 717.74                                                    & 667.83                                               &  & 41.14                                                                               & 41.86                                                                                            & 21.61                                                                                            & 660.40                                                    & 103.74                                               & 8.68                                                                                    \\
				& 60                        & 37.54                                                                               & 41.03                                                                                            & 31.96                                                                                            & 734.35                                                    & 1448.84                                              &  & 37.61                                                                               & 41.43                                                                                            & 20.43                                                                                            & 666.42                                                    & 442.19                                               & 10.19                                                                                   \\
				& 65                        & 35.93                                                                               & 39.40                                                                                            & 30.03                                                                                            & 755.87                                                    & 2269.15                                              &  & 35.69                                                                               & 40.41                                                                                            & 19.07                                                                                            & 680.04                                                    & 903.41                                               & 11.15                                                                                   \\
				& 70                        & 34.01                                                                               & 37.28                                                                                            & 28.82                                                                                            & 775.86                                                    & 2772.25                                              &  & 34.40                                                                               & 38.26                                                                                            & 17.39                                                                                            & 687.97                                                    & 2065.62                                              & 12.78                                                                                   \\
				& 75                        & 31.63                                                                               & 37.17                                                                                            & 29.86                                                                                            & 792.34                                                    & 3470.87                                              &  & 32.35                                                                               & 37.49                                                                                            & 15.23                                                                                            & 688.73                                                    & 2636.67                                              & 15.04                                                                                   \\ \cmidrule(lr){1-14} 
			\end{tabular}
		}
	\end{table}

	The results for the two-stage stochastic programming model and the dynamic rerouting approach are reported as percentages of the total value of assets impacted. The objective value ("Total value of assets protected") is the summation of the actual values of protected assets at stage one and the expected values of them at stage two.
	
	Table \ref{result} reveals that the dynamic rerouting approach usually protects more value in the first stage and the scenario with highest probability of occurrence (scenario 1) compared to the two-stage stochastic model. Planning protective tasks based on the first stage and scenario 1 means losing potential value of assets in scenario 2. However, the two-stage stochastic programming takes advantage of integrating all scenarios and maximises the total expected value of the serviced assets. Therefore, the two-stage stochastic programming may protect less value at the first stage and scenario 1 compared to the rerouting approach, but the total expected value of serviced assets in scenario 2 compensates for the loss.
	
	As is to be expected for both methods, it is seen in Table \ref{result} that computational time increases with increasing assets. Less obvious is that more time is required to solve instances with a larger number of vehicles. Dispatching more vehicles increases the complexity of the IMT's task but more assets are serviced. 
	The average computational time for set1 is around 340 second compared to 1416 seconds for set2. 
	In addition to the computational times the discrepancy in the objective values of the two approaches increases as more assets are involved. Increasing the number of resources in set2 and the problem size from 50 to 75 widens the objective function gap between the two methods from 5.91\% to 15.04\%. 
	Overall, the results indicates that the two-stage stochastic approach can improve the quality of the solution but that it comes at a computational cost. For less than 60 assets this computational cost is not a concern as solution times are still satisfactory for operational purposes. For cases involving a larger number of assets the rerouting approach would need to be employed despite the reduction in solution quality. However, that might not always be the case as different IMTs are responsible for various regions. This will cause limited number of key assets being assigned to each IMT for protection.
	
	\section{Extension} \label{extension}
	
	In section \ref{sec:formulation} a new two-stage stochastic approach proposed to solve the problems with uncertain time of changes. However, this would be more complicated when more than two scenarios are involved. 
	Assuming a problem comprised of four scenarios in which changes may take place at any of the following times: 2:00pm (first scenario), 2:30pm (second scenario), 3:00pm (third scenario) or 3:30pm (forth scenario). 
	Regardless of the time of change all scenarios share the same decision variables at 2:00pm (first stage variables). 
	The complexity arises after the staging time when the second, third and forth scenarios should share the same values of variables at some of the times.
	For instance, all three scenarios (2, 3 and 4) should take the same values for variables associated with the time between 2:00pm and 2:30pm. 
	This is owing to mutual visits after the staging time and by the second scenario occurrence	that would impact on all subsequent scenarios.
	The same condition applies to the variables related to the time between 2:30pm and 3:00pm that have to take the same values for scenarios occurring afterwards (scenarios 3 and 4).
	Our model can be extended to handle derived complexities in such situations, to consider a larger number of scenarios.
	
	Assuming that there are $F$ possible scenarios ordered based on the their time of occurrence in set $\Xi$ as:
	\begin{equation}\label{25}	 
	\Xi := \{\xi_1,\xi_2,\dots,\xi_F \};
	\end{equation}
	We define $TO_c$ as the time of occurrence for scenario $c$, therefore:
	\begin{equation}\label{26}	 
	TO_1 < TO_2 < \dots < TO_F. 
	\end{equation}
	The following constrains should be implemented to enforce equality of shared decision variables. Noting that $ST$ in constraints (\ref{6}) and (\ref{7}) needs to be replaced by $TO_1$. 
	\begin{equation}\label{27}	 
	TO_c > S_i^{s}(\xi_c)-M(1-\Gamma_i(\xi_c)) , \qquad \forall i\in N,  \xi_c \in \Xi \setminus{\xi_1};
	\end{equation}
	\begin{equation}\label{28}	 
	TO_c < S_i^{s}(\xi_c)+M\times\Gamma_i(\xi_c) , \qquad \forall i\in N,  \xi_c \in \Xi \setminus{\xi_1}. 
	\end{equation}
	Note that the variable $\Gamma_i(\xi_c)$ is equal to 1 if node $i$ under scenario $\xi_c$ is visited before the time $TO_c$. Using $\Gamma_i(\xi_c)$ we can identify nodes that are visited in shared time intervals. Writing below constraints for shared variables would allow us to assign the same values to them. $Y_i^{s}(\xi)$ and $z_{ijq}^{s}(\xi)$ are linked to $S_i^{s}(\xi)$ and $X_{jiq}^s(\xi)$, therefore we do not need to add constraints to assign same values to them in the shared time of scenarios.
	\begin{equation}\label{29}	 
	S_i^{s}(\xi_c) \le S_i^{s}(\xi_{c^\prime})+M(1-\Gamma_i(\xi_c)) , \qquad \forall i\in N,  \xi_c,\xi_{c^\prime}  \in \Xi \setminus{\xi_1}, c^\prime > c 
	\end{equation}
	\begin{equation}\label{30}	 
	S_i^{s}(\xi_c) \ge S_i^{s}(\xi_{c^\prime})+M(1-\Gamma_i(\xi_c)) , \qquad \forall i\in N,  \xi_c,\xi_{c^\prime}  \in \Xi \setminus{\xi_1}, c^\prime > c
	\end{equation}
	\begin{equation}\label{31}	
		X_{jiq}^s(\xi_c) \le X_{jiq}^s(\xi_{c^\prime})+M(1-\Gamma_i(\xi_c)) , \qquad \forall i,j\in N, q \in Q, \xi_c,\xi_{c^\prime}  \in \Xi \setminus{\xi_1}, c^\prime > c
		\end{equation}
		\begin{equation}\label{32}
		X_{jiq}^s(\xi_c) \ge 	X_{jiq}^s(\xi_{c^\prime})+M(1-\Gamma_i(\xi_c)) , \qquad \forall i,j\in N, q \in Q, \xi_c,\xi_{c^\prime}  \in \Xi \setminus{\xi_1}, c^\prime > c
	\end{equation}
\begin{equation}\label{33}	 
\Gamma_i(\xi_c) \in \{0,1\}, \qquad \forall i\in N,  \xi_c, \in \Xi.
	\end{equation}
\section{Conclusion} \label{conclusion}
	Wildfires are often hazardous to key assets whose loss can disrupt community life for months and be costly to replace. During an uncontrollable wildfire IMT's deploy response vehicles to mitigate the risk of losing crucial assets. Planning, coordinating and managing protective activities using limited resources (i.e., personnel and equipment) in an optimal way is critically important. Any strategy would not be optimal without utilising all available knowledge including the forecast of a wind change even if there is some uncertainty about the timing of that change. This is an important problem that arises frequently with wildfires in south-eastern Australia.
	In this paper we developed an extension to the asset protection problem. A two-stage stochastic programming approach was developed to handle the unusual feature of uncertainty in the timing of changes in conditions. In our study the changed conditions refer to wind direction and velocity. These changes determine new time windows during which assets must be serviced and hence the optimal deployment schedule and routing of vehicles. To the best of our knowledge this is the first two-stage stochastic programming problem dealing with uncertainty in the timing of a change.
	
	In addition to a case study, the proposed approach was tested with extensive computational experiments. To evaluate the efficacy of our solution scheme, results were compared to the dynamic rerouting approach. As expected the two-stage stochastic program produced better solutions than the rerouting approach. This difference became more pronounced as the complexity of the problem increased with more assets and vehicles. The performance of the rerouting approach in terms of computational time was superior to the stochastic programming approach but this only becomes relevant when a very large number of assets are involved. Given the scale of responsibility of most IMT's this limit is not expected to be of any operational consequence.
	
	Although the discussion focused on asset protection and change in wind direction, other logistic operations can be accommodated by this modelling approach. Along with weather change, other disruptions with uncertain time of occurrence include problems such as network accessibility, and change in travel time. These can be incorporated into the model in  a similar way.  
	
	It is believed that the proposed approach will be useful for other logistic problems where there is uncertainty about the timing of an event that is known will happen. This frequently occurs, for example, with picking up a passenger at an airport. It is known that the passenger will arrive and needs to be picked up but the timing of the pick-up is not precise due to uncertainties over the time they will need to collect baggage and pass through customs and immigration. This might have implications for scheduling of shuttle busses or fleets of taxis. Project managers schedule activities on multiple construction sites but the completion time of some stages are more dependent on weather and other conditions than others. The completion will happen but its timing is uncertain. Our approach might be a better way to deal with this uncertainty than current approaches.

	\appendix
	
	\section{A MILP model for deterministic asset protection problem}\label{appendixA1}
	
	The asset protection problem is formulated as the following mixed integer programming model (see \cite{roozbeh2018adaptive}).
	
	\begin{equation}\label{1a}
	Max \sum_{i\in N} \nu_iY_i 
	\end{equation}
	\begin{equation}\label{2a}
	s.t.:\sum_{j \in \delta^+_q(0)} X_{0jq} = \sum_{i \in \delta^-_q(N+1)} X_{i(N+1)q}, \quad q\in Q,
	\end{equation}
	\begin{equation}\label{3a}
	\sum_{i \in \delta^-_q(k)} X_{ikq} = \sum_{j \in \delta^+_q(k)} X_{kjq}, \; k=1,\dots,N, \quad q\in Q,
	\end{equation}
	\begin{equation}\label{4a}
	r_{kq}Y_{k} = \sum_{i \in \delta^-_q(k)} X_{ikq}, \; k=1,\dots,N, \quad q\in Q,
	\end{equation}
	\begin{equation}\label{5a}
	X_{ijq}\leq \kappa_{q} z_{ijq}, \;(i,j)\in A, \quad q\in Q,
	\end{equation}
	\begin{equation}\label{6a}
	S_i+t_{ijq}+a_i-S_j\leq M(1-z_{ijq}), \;(i,j)\in A,  \; q\in Q,
	\end{equation}
	\begin{equation}\label{7a}
	o_i \leq S_i, \; i=1,\dots,N,
	\end{equation}
	\begin{equation}\label{8a}
	S_i \leq c_i, \; i=1,\dots,N,
	\end{equation}
	\begin{equation}\label{9a}
	x_{ijq} \in \{0,1,\dots,p_q\}, \; (i,j)\in A,
	\end{equation}
	\begin{equation}\label{10a}
	y_i,z_{ijq}\in \{0,1\}, \; (i,j)\in A.
	\end{equation}
	
	The objective function maximises the total value of the protected assets. Constraints (\ref{2a} and \ref{3a}) ensure flow conservation for depot and all other assets. Constraint (\ref{4a}) guarantees protection of an asset upon satisfaction of its protection requirements. Constraint (\ref{5a}) defines an upper bound for different types of resources. Constraint (\ref{6a}) ensures that an asset can be visited when the protection requirements of the previous location has been satisfied and there is enough time to reach the next asset. Equations (\ref{7a}) and (\ref{8a}) ensure that the time window constraints are met. Integer and binary conditions are defined in constraints (\ref{9a}) and (\ref{10a}).

	\section{Supplementary data for the Black Saturday case study}\label{appendixB}
	A weather change is forecast after the initial trip scheduled for the Black Saturday case study (section \ref{case}). In table \ref{assetlist},  assets impacted as a result of the change in the weather condition are identified on the map. The protection requirement vector shows the number of resources of each type needed for accomplishing the protection task. Protection requirement vectors are assigned randomly and time windows are defined approximately based on The Victorian Bushfires Royal Commission report (\cite{whittaker2009victorian}). Some assets are impacted in multiple scenarios and others in one scenario, while those that are not at risk require no protection and have the impact time of 0.
	\newpage
	\begin{table}[h]
		\centering
		\caption{List of assets being impacted in section \ref{case}. For sensitivity reasons asset names are not mentioned. }
		\label{assetlist}
		\resizebox{1\columnwidth}{!}{
			\begin{tabular}{@{}lccccc@{}}
				\toprule
				\# & Asset ID                           & \begin{tabular}[c]{@{}c@{}}Protection \\ requirement vector\end{tabular} & \begin{tabular}[c]{@{}c@{}}Time of being impacted \\ at stage 1\end{tabular} & \begin{tabular}[c]{@{}c@{}}Time of being impacted \\ at scenario 1 (stage 2)\end{tabular} & \begin{tabular}[c]{@{}c@{}}Time of being impacted \\ at scenario 2 (stage 2)\end{tabular} \\ \midrule
				1  & A1                   & \textless{}2, 1, 0\textgreater{}                                         & 0                                                                            & 0                                                                                         & 7                                                                                         \\
				2  & A2                  & \textless{}2, 0, 1\textgreater{}                                         & 3                                                                            & 0                                                                                         & 0                                                                                         \\
				3  & A3                       & \textless{}1, 0, 2\textgreater{}                                         & 0                                                                            & 7.5                                                                                       & 7.5                                                                                       \\
				4  & A4                 & \textless{}0, 2, 1\textgreater{}                                         & 4                                                                            & 0                                                                                         & 0                                                                                         \\
				5  & A5        & \textless{}1, 0, 2\textgreater{}                                         & 0                                                                            & 10                                                                                        & 0                                                                                         \\
				6  & A6            & \textless{}1, 1, 1\textgreater{}                                         & 0                                                                            & 7                                                                                         & 7                                                                                         \\
				7  & A7                    & \textless{}2, 1, 0\textgreater{}                                         & 4                                                                            & 0                                                                                         & 0                                                                                         \\
				8  & A8         & \textless{}2, 0, 1\textgreater{}                                         & 0                                                                            & 0                                                                                         & 8                                                                                         \\
				9  & A9            & \textless{}1, 2, 0\textgreater{}                                         & 0                                                                            & 11.5                                                                                      & 0                                                                                         \\
				10 & A10              & \textless{}1, 0, 2\textgreater{}                                         & 0                                                                            & 8                                                                                         & 8                                                                                         \\
				11 & A11                        & \textless{}1, 0, 2\textgreater{}                                         & 5                                                                            & 0                                                                                         & 0                                                                                         \\
				12 & A12      & \textless{}1, 2, 0\textgreater{}                                         & 5.5                                                                          & 0                                                                                         & 0                                                                                         \\
				13 & A13         & \textless{}1, 0, 2\textgreater{}                                         & 6                                                                            & 0                                                                                         & 0                                                                                         \\
				14 & A14             & \textless{}0, 2, 1\textgreater{}                                         & 0                                                                            & 0                                                                                         & 7                                                                                         \\
				15 & A15                     & \textless{}2, 1, 0\textgreater{}                                         & 0                                                                            & 8                                                                                         & 8                                                                                         \\
				16 & A16                      & \textless{}2, 0, 1\textgreater{}                                         & 6.5                                                                          & 0                                                                                         & 0                                                                                         \\
				17 & A17            & \textless{}2, 0, 1\textgreater{}                                         & 4                                                                            & 0                                                                                         & 0                                                                                         \\
				18 & A18                  & \textless{}2, 0, 1\textgreater{}                                         & 0                                                                            & 7                                                                                         & 0                                                                                         \\
				19 & A19                  & \textless{}2, 0, 1\textgreater{}                                         & 0                                                                            & 12                                                                                        & 0                                                                                         \\
				20 & A20 & \textless{}1, 2, 0\textgreater{}                                         & 0                                                                            & 0                                                                                         & 8                                                                                         \\
				21 & A21                      & \textless{}1, 2, 1\textgreater{}                                         & 0                                                                            & 10                                                                                        & 10                                                                                        \\
				22 & A22                      & \textless{}1, 2, 1\textgreater{}                                         & 0                                                                            & 11                                                                                        & 0                                                                                         \\
				23 & A23                    & \textless{}1, 2, 1\textgreater{}                                         & 0                                                                            & 11                                                                                        & 11                                                                                        \\
				24 & A24              & \textless{}1, 2, 1\textgreater{}                                         & 0                                                                            & 11.5                                                                                      & 0                                                                                         \\
				25 & A25                 & \textless{}2, 0, 1\textgreater{}                                         & 0                                                                            & 9.5                                                                                       & 0                                                                                         \\ \bottomrule
			\end{tabular}
		}
	\end{table}

	
	

	{
		
	}

	\newpage

	\bibliographystyle{elsarticle-harv}
	\bibliography{p3}
	
\end{document}